\documentclass[12pt]{article}

\usepackage{epsfig}
\usepackage{amssymb}
\usepackage{amsmath}
\usepackage{amsthm}
\usepackage{latexsym}

\def\P{{\mathbb P}}
\newtheorem{theorem}{Theorem}[section]
\newtheorem{lemma}[theorem]{Lemma}
\newtheorem{corollary}[theorem]{Corollary}

\title{
Large-N Limit of Crossing Probabilities, Discontinuity,
and Asymptotic Behavior \\ of Threshold Values in Mandelbrot's\\
Fractal Percolation Process}

\author{
{Erik I. Broman}
\thanks{Department of Mathematics, Chalmers University of Technology,
S-412 96 G\"oteborg, Sweden. E-mail: broman\,@\,math.chalmers.se}
\thanks{The work of this author was carried out while at the Department
of Mathematics of the Vrije Universiteit Amsterdam.} \,\,\,\,\,\,\,\,\,\,
{Federico Camia}
\thanks{Department of Mathematics, Vrije Universiteit Amsterdam,
De Boelelaan 1081a, 1081 HV Amsterdam, The Netherlands. E-mail: fede\,@\,few.vu.nl}
\thanks{Partially supported by a VENI grant of the NWO.}\\
}

\date{}

\begin{document}

\maketitle

\begin{abstract}

We study Mandelbrot's percolation process in dimension $d \geq 2$.
The process generates random fractal sets by an iterative procedure
which starts by dividing the unit cube $[0,1]^d$ in $N^d$ subcubes,
and independently retaining or discarding each subcube with
probability $p$ or $1-p$ respectively. This step is then repeated
within the retained subcubes at all scales. As $p$ is varied,
there is a percolation phase transition in terms of paths for all
$d \geq 2$, and in terms of $(d-1)$-dimensional ``sheets" for all
$d \geq 3$.


For any $d \geq 2$, we consider the random fractal set produced
at the path-percolation critical value $p_c(N,d)$, and show that
the probability that it contains a path connecting two opposite
faces of the cube $[0,1]^d$ tends to one as $N \to \infty$.
As an immediate consequence, we obtain that the above probability
has a discontinuity, as a function of $p$, at $p_c(N,d)$ for all
$N$ sufficiently large. This had previously been proved only for
$d=2$ (for any $N \geq 2$). For $d \geq 3$, we prove analogous
results for sheet-percolation.

%

In dimension two, Chayes and Chayes proved that $p_c(N,2)$
converges, as $N \to \infty$, to the critical density $p_c$ of site
percolation on the square lattice. Assuming the existence of the
correlation length exponent $\nu$ for site percolation on the square
lattice, we establish the speed of convergence up to a logarithmic
factor. In particular, our results imply that
$p_c(N,2)-p_c=(\frac{1}{N})^{1/\nu+o(1)}$ as $N \to \infty$, showing
an interesting relation with near-critical percolation.

\end{abstract}



\section{Introduction and First Results}

In this paper we are concerned with a continuum percolation model,
first introduced in~\cite{Mandelbrot}, which is known as Mandelbrot's
fractal percolation process and defined as follows.
For any integers $d \geq 2$ and $N \geq 2,$ we start by partitioning
the unit cube $[0,1]^d \subset {\mathbb R}^d$ into $N^d$ subcubes
of equal size. Given $p \in [0,1]$ and a subcube, we then retain the
subcube with probability $p$ and discard it with probability $1-p$.
This is done independently for every subcube of the partition. This
gives us a random set $C_N^1=C_N^1(d,p) \subset [0,1]^d.$
Consider any retained (assuming that $C_N^1 \neq \emptyset$) subcube
$B$ in $C_N^1.$ We can repeat the described procedure on a smaller
scale by partitioning $B$ into $N^d$ further subcubes, discarding or
retaining them as above. We do this for every retained subcube of
$C_N^1.$ This will yield a new random set $C_N^2 \subset C_N^1.$
Iterating the procedure on every smaller scale yields an infinite
sequence of random sets $\ldots \subset C_N^{k+1} \subset C_N^k
\subset \ldots \subset [0,1]^d.$ It is easy to see that
$C_N:=\cap_{k=1}^{\infty} C_N^k$ is a well defined random object.

We denote by $CR(C_N)$ the event that $C_N$ contains a connected component
intersecting both the ``left-hand face" $\{ 0 \} \times [0,1]^{d-1}$ and
the ``right-hand face" $\{ 1 \} \times [0,1]^{d-1}$ of $[0,1]^d$ (we will
call such a connected component a path or crossing). We can then define
\begin{equation} \label{eqn10}
\theta_{N,d}(p):={\mathbb P}_p(CR(C_N)).
\end{equation}
We could of course have used any two opposite faces in this definition.
Furthermore, we define
\[
p_c(N,d):=\inf\{p:\theta_{N,d}(p)>0\}.
\]
By a standard coupling argument it is easy to see that
$\theta_{N,d}(p)$ is non-decreasing in $p$ for every $N$.

Several authors studied various aspects of Mandelbrot's fractal
percolation, including the Hausdorff dimension of the limiting
fractal set, as detailed in~\cite{DG}, and the possible
existence of paths~\cite{CCD,DM,Meester,CCGS} and $(d-1)$-dimensional
``sheets" \cite{CCGS,Orzechowski,MPV} traversing the unit cube between
opposite faces. Dekking and Meester~\cite{DM} proposed a ``morphology"
of random Cantor sets comprising several ``phases" in which a set
can be.

Chayes and Chayes~\cite{ChCh} considered the behavior of $p_c(N,2)$
for large $N$, and proved that it converges to the critical value for
Bernoulli site percolation on the square lattice, as $N \to \infty$.
This was generalized by Falconer and Grimmett~\cite{FG1,FG2} to all
dimensions $d \geq 3$, where the limiting value however is not the
critical value of site percolation on the (hyper-)cubic lattice, but
that of site percolation on a related lattice (see Section~\ref{sec-enh-dim}).
Orzechowski~\cite{Orzechowski} proved an analogous results for
``sheet" percolation.

Falconer and Grimmett~\cite{FG1,FG2} obtained their result as a
consequence of proving that, for all $d\geq2$,
$\lim_{N \to \infty}\theta_{N,d}(p)=1$ for all $p$ strictly larger
than the relevant (depending on the dimension) critical value.
An interesting question, addressed in Theorem~\ref{thm3} below,
concerns the large-$N$ limit of the path-crossing probability
\emph{at} $p_c(N,d)$. Let us try to motivate that question.

For the sake of this discussion, let us consider the model in two
dimensions, where it is known~\cite{CCD,DM} that $\theta_{N,2}(p)$
is discontinuous at $p_c(N,2)$ for all $N$. Analogously, if we tile
the plane with independent copies of the system, the probability
that there exists an infinite path intersecting the unit square
$[0,1]^2$ has a discontinuity at $p_c(N,2)$ for all $N$~\cite{CCD}.
Such behavior prompted the authors of~\cite{CCD} to call the phase
transition in Mandelbrot's fractal percolation first order or
discontinuous (although one of the authors later appeared to change his
mind on the issue -- see p.~130 of~\cite{lchayes}).
Regardless of the label, the discontinuity of $\theta_{N,2}(p)$ (and, more
generally, of $\theta_{N,d}(p)$ for $d\geq2$) is an important feature of
the phase transition, and it is natural to ask what happens to it as
$N \to \infty$.


We point out that the fact that
$\lim_{N \to \infty}\theta_{N,d}(p)=1$ for all $p$ strictly larger than the
relevant (depending on the dimension) critical value does not shed much light
on the large-$N$ limit of the path-crossing probability \emph{at} the critical
point. Consider, for instance, the crossing probability of a square in (ordinary)
two-dimensional Bernoulli percolation in the scaling limit (i.e., as the lattice
spacing goes to zero) or equivalently in the limit of larger and larger squares.
Such a limit is needed to ``see" the phase transition and produce a situation
that can be compared to the one in Mandelbrot's fractal percolation. However, in
Bernoulli percolation, in the scaling limit or in the limit of larger and larger
squares, the crossing probability of a square converges to one above the critical
point, to zero below it, and to a value strictly between zero and one \emph{at} the
critical point, showing a behavior that differs from the one of the next theorem.

%

\begin{theorem} \label{thm3}
We have that for any $d \geq 2,$
\[
\lim_{N \rightarrow \infty}\theta_{N,d}(p_c(N,d))=1.
\]
\end{theorem}

When proving Theorem~\ref{thm3}, in Section~\ref{sec-proof-thm3}, we will
use a result from~\cite{AG} concerning ``enhancement" percolation (see
also~\cite{G}). We will in fact need a slightly modified version of the result;
in particular, we will be looking at crossing probabilities in a ``diminishment"
percolation model (to be defined precisely in Section~\ref{sec-enh-dim}).

The following corollary is an immediate consequence of the previous
theorem. For $d=2$ and all $N\geq2$, an elegant and rather elementary
proof of this fact can be found in~\cite{DM}.

\begin{corollary} \label{discontinuity}
For any $d \geq 2,$ there exists an $N_0=N_0(d)$ such that the
function $\theta_{N,d}(p)$ is discontinuous at $p_c(N,d)$ for
$N \geq N_0.$
\end{corollary}

\noindent{\bf Remark.} For $d=3$ and all $N\geq2$, a similar result had
been proved in~\cite{CCGS}, but with the probability of crossing a square
replaced by that of crossing a $2\times2\times1$ rectangle between the
$2\times2$ faces. While the result in~\cite{CCGS} is more satisfactory
because it applies to all $N\geq2$, it is limited to a three-dimensional
system with a special geometry. On the contrary our result is valid in
all dimensions and, although we use the unit cube for simplicity, the
choice of that particular geometry is not important for our arguments. \\


When $d \geq 3$, one can consider the existence of $(d-1)$-dimensional
``sheets'' crossing the unit cube $[0,1]^d$.
We denote by $SH(C_N)$ the event that $C_N$ contains a surface
separating the faces of $[0,1]^d$ perpendicular to the first coordinate
direction. We can then define
\begin{equation} \label{critical-sheet}
\tilde\theta_{N,d}(p):={\mathbb P}_p(SH(C_N)).
\end{equation}
We could of course have used any two opposite faces in this definition.
Furthermore, we define
\[
\tilde{p}_c(N,d):=\inf\{p:\tilde\theta_{N,d}(p)>0\}.
\]
By a standard coupling argument it is easy to see that
$\tilde\theta_{N,d}(p)$ is non-decreasing in $p$ for every $N$.

Our next result concerns the large-$N$ limit of the sheet-crossing probability
at the critical point, and is analogous to Theorem~\ref{thm3}.
\begin{theorem} \label{thm-sheet}
We have that for any $d \geq 3,$
\[
\lim_{N \rightarrow \infty}\tilde\theta_{N,d}(\tilde{p}_c(N,d))=1.
\]
\end{theorem}

When proving this theorem, in Section~\ref{sec-proof-thm-sheet}, we
will again use a slight modification of a result from~\cite{AG}
concerning enhancement percolation.

The following corollary is an immediate consequence of the previous
theorem.

\begin{corollary} \label{discontinuity-sheet}
For any $d \geq 3,$ there exists an $\tilde{N}_0=\tilde{N}_0(d)$ such that
the function $\tilde\theta_{N,d}(p)$ is discontinuous at $\tilde{p}_c(N,d)$
for $N \geq \tilde{N}_0.$
\end{corollary}


\section{Asymptotic Behavior of $p_c(N,2)$} \label{section2d}

Our last result concerns the asymptotic behavior of $p_c(N,2)$ as $N
\to \infty$. It is known~\cite{ChCh} that, for $d=2$, $p_c(N,2)$
tends to the critical density $p_c=p_c({\mathbb Z}^2,site)$ of
(Bernoulli) site percolation (see Section~\ref{sec-enh-dim} for the exact
definition) on the square lattice, as $N \to \infty$. Assuming the existence
of the correlation length exponent $\nu$ for two-dimensional site percolation
on the square lattice, we give bounds on the speed of convergence.


There are several equivalent ways to introduce the concept of correlation
length in percolation.
In this section we will assume that $p>p_c$, since this is the only case
we are interested in.
Following Kesten (see~\cite{Kesten},~\cite{ChN} and the references therein),
we define
\begin{equation}
L_{\delta}(p) := \min \{ N : \P_p(CR(C^1_N)) \geq 1-\delta  \}.
\end{equation}

The $\delta$ in the definition is unimportant, since for any $\delta, \delta' \in (0,1/2)$
we have~\cite{Kesten,Nolin} that the ration $L_{\delta}(p)/L_{\delta'}(p)$ is bounded away
from $0$ and $\infty$ as $p \to p_c$. In view of this, in what follows we will consider
$\delta$ to be a fixed number smaller than $1/2$. $L_{\delta}(p)$ is believed to behave
like $(p-p_c)^{-\nu}$ as $p \downarrow p_c$, with $\nu=4/3$ (see, e.g., \cite{Kesten} and
references therein). More precisely, taking $\delta<1/2$, we make the following assumption.\\

\noindent{\bf H1}: There are two constants $c_1=c_1(\delta)$ and $c_2=c_2(\delta)$
such that $c_1 (p-p_c)^{-4/3} \leq L_{\delta}(p) \leq c_2 (p-p_c)^{-4/3}$ as
$p \downarrow p_c$.\\

Based on that assumption, the following theorem is proved in Section~\ref{sec-proof-2d}.
\begin{theorem} \label{thm4} Let $d=2$, and assume that {\bf H1} holds.
For any $C<\infty$, there exists an $N_C<\infty$ such that for all $N
\geq N_C,$
\begin{equation} \label{eqn15}
p_c(N,2) \geq p_c+C\left( \frac{1}{N}\right)^{3/4}.
\end{equation}
Furthermore, there exist constants $D<\infty$ and $N_1<\infty$ such that,
for all $N \geq N_1$,
\begin{equation} \label{eqn16}
p_c(N,2) \leq p_c+D\left( \frac{\log N}{N}\right)^{3/4}.
\end{equation}
\end{theorem}


The limit $N \to \infty$ can be compared to the scaling limit $\varepsilon \to 0$
in ordinary percolation on a lattice with lattice spacing $\varepsilon$. The reason
is that it removes the intrinsic scale of the system due to the size of the largest
subdivisions in the same way as the limit $\varepsilon \to 0$ in ordinary percolation
removes the intrinsic scale of the system due to the lattice. Theorem~\ref{thm4}
and its proof show an interesting (and at first site maybe surprising) connection
between the large-$N$ limit of fractal percolation and the ``near-critical" scaling
limit of ordinary percolation. The latter corresponds to letting $\varepsilon \to 0$
in a percolation model with lattice spacing $\varepsilon$ and density of open sites
(or bonds) $p=p_c+\lambda \varepsilon^{3/4}$ for some fixed $\lambda$ (see, e.g.,
\cite{bcks,cfn1,cfn2,nw}).

To motivate the connection, roughly speaking, one can argue in the following way.
Take $N$ very large, and consider a system with retention probability equal to $p_c$.
After the first iteration the system looks critical, meaning that there is positive
probability to see left-to-right crossings, but that such crossings are tenuous
and are not able to sustain the remaining iterations. This is not surprising since
we know that $p_c(N,2)>p_c$ for all $N$. On the other hand, if $p>p_c$, for $N$
sufficiently large, a system with retention probability $p$ looks very supercritical
after the first iteration (in other words, $N \gg L_{\delta}(p)$). Chayes and
Chayes~\cite{ChCh} used this fact to show that the system has enough connections
to sustain all further iterations. As a consequence, one obtains that $p_c(N,2)$
becomes arbitrarily close to $p_c$ as $N \to \infty$. But what is the speed of
convergence? As $N \to \infty$, $p_c(N,2)$ must be getting close to $p_c$ at the
right speed: not too fast, or else after the first iteration the system would look
``too critical," and not too slowly, or else it would look ``too supercritical."
As Theorem~\ref{thm4} and its proof show, the right speed is related to near-critical
percolation. \\

\noindent{\bf Remark}.
Although {\bf H1} is widely believed to hold (at least) for site and bond
percolation on regular lattices, it has so far not been proved for any
percolation model. However, a weaker result has been proved~\cite{SW} for
site percolation on the triangular lattice. With the notation introduced
above, the latter can be stated as follows.\\

\noindent{\bf H2}: $L_{\delta}(p) = (p-p_c)^{-4/3 + o(1)}$ as
$p \downarrow p_c$.\\

\noindent Assuming {\bf H2} instead of {\bf H1}, our results can be easily
adapted to obtain $p_c(N,2)-p_c=(\frac{1}{N})^{3/4+o(1)}$.\\

\noindent{\bf Remark}.
The definition of correlation length we use in this paper is equivalent
to other standard definitions (see, e.g., \cite{Kesten,ChN,Nolin}), including
the one that comes from considering the probability $\tau^f_p(x)$ that the
origin and $x$ belong to the same finite $p$-open cluster.
The correlation length $\xi(p)$ related to $\tau^f_p(x)$ is defined by
\begin{equation} \label{correlation-length} \nonumber
\xi(p)^{-1} := \lim_{|x| \to \infty} \left\{ - \frac{1}{|x|} \log
\tau^f_p(x) \right\}.
\end{equation}
In particular, for every $\delta>0$, there are two constants, $c'_1=c'_1(\delta)$
and $c'_2=c'_2(\delta)$, such that $c'_1 L_{\delta}(p) \leq \xi(p) \leq c'_2 L_{\delta}(p)$
for all $p>p_c$ (see, e.g., the Appendix of~\cite{ChN}). \\

\noindent{\bf Remark:} In light of Theorem \ref{thm4}, we find the concluding remarks
of~\cite{ChCh} unclear. On p.~L505 of~\cite{ChCh}, the authors claim that: ``In the
large-N regime, our proof demonstrates that if $N$ is a huge (but fixed) multiple of
the correlation length of the density-p site model, then $p>p_c(N)$." However, if one
accepts hypothesis {\bf H1}, this would imply that, for $N$ large enough,
$p_c(N,2) < p_c + \frac{const}{N^{3/4}}$, which contradicts the first statement of
Theorem~\ref{thm4}. If one instead accepts {\bf H2}, then it is not clear what a
multiple of the correlation length means.

\section{Enhancement and Diminishment Percolation} \label{sec-enh-dim}

In this section we introduce and discuss a particular example of enhancement
percolation~\cite{AG}, and a closely related example of ``diminishment"
percolation (see the comment after Theorem~3.16 on p.~65 of~\cite{G}).
Both models will be used later in the proofs of the main results, where
they will be compared (or more precisely, coupled) to the fractal
percolation model.


Before we can proceed we need to introduce the two graphs on which the
enhancement and diminishment percolation models will be defined. For
$d \geq 2$, let ${\mathbb L}^d$ be the $d$-dimensional lattice with
vertex set ${\mathbb Z}^d$ and with edge set given by the adjacency
relation: $(x_1, \ldots, x_d)=x \sim y=(y_1, \ldots, y_d)$
if and only if $x \neq y$, $|x_i-y_i| \leq 1$ for all $i$
and $x_i=y_i$ for at least one value of $i$. Let ${\mathbb M}^d$
be the $d$-dimensional lattice with vertex set ${\mathbb Z}^d$ and
with edge set given by the adjacency relation:
$(x_1, \ldots, x_d)=x \sim y=(y_1, \ldots, y_d)$ if and only if
$x \neq y$ and $|x_i-y_i| \leq 1$ for all $i$.
For $d=2$, ${\mathbb L}^2$ coincides with the square lattice and
${\mathbb M}^2$ with its close-packed version, obtained by adding
diagonal edges inside the faces of the square lattice. For $d \geq 3$,
${\mathbb M}^d$ contains ${\mathbb L}^d$ as a strict sublattice, which
in turns contains the (hyper-)cubic lattice as a strict sublattice.


We consider Bernoulli site percolation on ${\mathbb L}^d$ and
${\mathbb M}^d$ with parameter $p \in [0,1]$: each site is declared open
(represented by a 1) with probability $p$ and closed (represented by a 0)
with probability $1-p$, independently of all other sites.
For any $N \geq 2,$ consider the box $\Lambda^d_N := \{1,\ldots,N\}^d
\subset {\mathbb Z}^d$.

For $d \geq 2$, we denote by $CR(\Lambda^d_N,{\mathbb L}^d)$ the event that there
is an open path in ${\mathbb L}^d$ crossing $\Lambda^d_N$ in the first coordinate
direction. An open path in ${\mathbb L}^d$ is a sequence of open sites such that
any two consecutive elements of the sequence are adjacent in ${\mathbb L}^d$. If
$CR(\Lambda^d_N,{\mathbb L}^d)$ occurs, we say that there is an open crossing of
$\Lambda^d_N$. Analogous definitions hold for $d \geq 3$ with ${\mathbb L}^d$
replaced by ${\mathbb M}^d$.

For $d \geq 2$, we define
\begin{equation} \label{phi}
\varphi^d_N(p):={\mathbb P}_p(CR(\Lambda^d_N,{\mathbb L}^d))
\end{equation}
and
\[
p_c({\mathbb L}^d):=\inf\{p: \lim_{N \rightarrow \infty }
\varphi^d_N(p)>0\}.
\]
As before, we could of course have considered crossings in any direction.



Analogously, for $d \geq 3$, we define
\begin{equation} \label{psi}
\psi^d_N(p):={\mathbb P}_p(CR(\Lambda^d_N,{\mathbb M}^d))
\end{equation}
and
\[
p_c({\mathbb M}^d):=\inf\{p: \lim_{N \rightarrow \infty }
\psi^d_N(p)>0\}.
\]


In~\cite{FG2} it is proved, for all $d \geq 2$, that $p_c({\mathbb
L}^d) \leq p_c(N,d)$ for all $N$, and that $$\lim_{N \rightarrow
\infty}p_c(N,d)=p_c({\mathbb L}^d).$$ (The same conclusions had been
reached earlier in~\cite{ChCh} for $d=2$.) Similar results are obtained
for the sheet percolation problem in dimension $d \geq 3$
in~\cite{Orzechowski}, where it is shown that $\tilde{p}_c(N,d) \geq
1 - p_c({\mathbb M}^d)$ for all $N$, and that
$$\lim_{N \rightarrow \infty}\tilde{p}_c(N,d) = 1 -
p_c({\mathbb M}^d).$$

In fact, the above results were stated for a different definition
of $p_c({\mathbb M}^d)$ and $p_c({\mathbb L}^d)$, defined using the
probability of the existence
of an infinite cluster. However, it is a standard result from
percolation theory (see, e.g., \cite{G} and references therein)
that the latter definition is equivalent to the one given above
using limits of crossing probabilities.

We are now ready to define our (stochastic) enhancement percolation
model with density $s$ of enhancement. We start by performing Bernoulli
site percolation with density $p\in[0,1]$ on ${\mathbb M}^d.$  For a site
$x=(x_1,\ldots, x_d)\in{\mathbb M}^d$, if $(x_1 + 1,x_2,\ldots, x_d)$
and $(x_1 - 1,x_2,\ldots, x_d)$ are both open, we make $x$ open with
probability $s$, regardless of its state before the enhancement. Note
that this has the effect of ``enhancing"the percolation configuration
by changing the state of some sites from closed to open. We denote by
$\P^{enh}_{p,s}$ the resulting probability measure.

It is immediate to see that our enhancement is essential in the
language of~\cite{AG} (i.e., there exist configurations that do not
have doubly-infinite open paths but such that a doubly-infinite open
path [the union of two disjoint, infinite open paths starting at
neighboring sites] is generated if the enhancement is activated at
the origin).

For $s>0$ fixed, we define
\[
\psi^d_N(p,s):={\mathbb P}^{enh}_{p,s}(CR(\Lambda^d_N,{\mathbb M}^d)).
\]
Note that the enhancement breaks the symmetry of the original percolation
model, and recall that $CR(\cdot)$ is the event that there is a crossing
in the first coordinate direction. In contrast to similar definitions above,
choosing the first coordinate direction here is important. Note also that
$\psi^d_N(p,0)=\psi^d_N(p)$, where the right hand side was defined in~(\ref{psi}).

The next lemma follows easily from the proof of Theorem~1 of~\cite{AG}.

\begin{lemma} \label{lemma-enhancement}
For fixed $s>0,$ there exists $\delta_1=\delta_1(s)>0$ so that,
for every $p$ such that $0 < p_c-p \leq \delta_1$,
\[
\lim_{N \rightarrow \infty}\psi^d_N(p,s) = 1.
\]
\end{lemma}

The proof of Theorem~1 of~\cite{AG} is based on a differential
inequality. In order to use the same strategy to prove
Lemma~\ref{lemma-enhancement}, we would need the following
differential inequality:
\begin{equation} \label{diff-ineq1}
\frac{\partial}{\partial s}\psi^d_N(p,s) \geq f_d(p,s)
\frac{\partial}{\partial p}\psi^d_N(p,s)
\end{equation}
with $f_d(p,s)$ continuous and strictly positive on $(0,1)^2$. This
can be easily obtained by applying (a straightforward generalization
of) Russo's formula (see \cite{AG} and \cite{G}) and the proof of
Lemma~2 of~\cite{AG} to the event $CR(\Lambda^d_N,{\mathbb M}^d)$.
Lemma~2 of~\cite{AG} concerns a different event, namely the
existence of an open path connecting the origin to the boundary of
the cube of side length $2N$ (say). However, the key ingredient of
the proof of the lemma is a local rearrangement of a percolation
configuration that involves changing the state of at most a bounded
number of sites, and is therefore not sensitive to the ``global
geometry" of the event.

We do not give a proof of~(\ref{diff-ineq1}) here, since it would
look almost identical to the proof of~(1.7) of~\cite{AG}, but for
the reader's convenience, we explain how to obtain
Lemma~\ref{lemma-enhancement} from~(\ref{diff-ineq1}). The latter
shows that $\psi^d_N(p,s)$ is a nonincreasing function of $t$ when
$(p,s) \equiv (p(t),s(t))$ satisfies
\begin{equation} \nonumber
\frac{d}{dt}(p,s) = (f_d(p,s),-1).
\end{equation}
In particular, since $f_d$ is strictly positive, for any $s_0>0$ there
is some $p_0<p_c({\mathbb M}^d)$ such that the forward orbit from
$(p_0,s_0)$ crosses the segment $\{ p=p_c({\mathbb M}^d), s \in [0,1]
\}$ and reaches some point $(p',s')$ with $p'>p_c({\mathbb M}^d)$.
Then, using the monotonicity in $s$, $\psi^d_N(p_0,s_0) \geq
\psi^d_N(p',s') \geq \psi^d_N(p',0) = \psi^d_N(p')$, which implies
the lemma.

Next, we define our (stochastic) diminishment percolation model with
density $s$ of diminishment. We start by performing Bernoulli site
percolation with density $p\in[0,1]$ on ${\mathbb L}^d$. For a site
$x=(x_1,\ldots, x_d)\in{\mathbb L}^d$, if all neighbors of $x$, except
possibly the ones with coordinates $(x_1 \pm 1,x_2,\ldots, x_d)$, are
closed, we make $x$ closed with probability $s$, regardless of its state
before the enhancement. Note that this has the effect of ``diminishing"
the percolation configuration by changing the state of some sites from
open to closed. We denote by $\P^{dim}_{p,s}$ the resulting probability
measure.

It is immediate to see that our diminishment is essential in the sense
that there exist configurations that have a doubly-infinite open path but
such that a doubly-infinite open path is not present after the enhancement
is activated at the origin.

For $s>0$ fixed, we define
\[
\varphi^d_N(p,s):={\mathbb P}^{dim}_{p,s}(CR(\Lambda^d_N,{\mathbb L}^d)).
\]
As above note that the diminishment breaks the symmetry of the original
percolation model so that again it matters that we consider crossings in
the first coordinate direction. Note also that
$\varphi^d_N(p,0)=\varphi^d_N(p)$, where the right hand side was defined
in~(\ref{phi}).

The next lemma can again be proved using arguments from~\cite{AG}.

\begin{lemma} \label{lemma-diminishment}
For fixed $s>0,$ there exists $\delta_2=\delta_2(s)>0$ so that for every $p$
such that $0 < p-p_c \leq \delta_2$,
\[
\lim_{N \rightarrow \infty}\varphi^d_N(p,s) =0.
\]
\end{lemma}

The lemma can be proved using the following differential inequality:

\begin{equation} \label{diff-ineq2}
\frac{\partial}{\partial s}\varphi^d_N(p,s) \leq - g_d(p,s)
\frac{\partial}{\partial p}\varphi^d_N(p,s),
\end{equation}
with $g_d(p,s)$ continuous and strictly positive on $(0,1)^2$.

Once again, we do not give a proof of~(\ref{diff-ineq2}), since it
would look almost identical to the proof of~(1.7) of~\cite{AG}.
To obtain Lemma~\ref{lemma-diminishment} from~(\ref{diff-ineq2}),
we can use the same strategy as before.
From~(\ref{diff-ineq2}) we see that $\varphi^d_N(p,s)$ is a nondecreasing
function of $t$ when $(p,s) \equiv (p(t),s(t))$ satisfies
\begin{equation} \nonumber
\frac{d}{dt}(p,s) = (-g_d(p,s),-1).
\end{equation}
In particular, since $g_d$ is strictly positive, for any $s_0>0$ there
is some $p_0>p_c({\mathbb L}^d)$ such that the forward orbit from
$(p_0,s_0)$ crosses the segment $\{ p=p_c({\mathbb L}^d), s \in [0,1]
\}$ and reaches some point $(p',s')$ with $p'<p_c({\mathbb L}^d)$.
Then, due to the monotonicity in $s$, $\varphi^d_N(p_0,s_0) \leq
\varphi^d_N(p',s') \leq \varphi^d_N(p',0) = \varphi^d_N(p')$,
which implies the lemma.

\section{Proofs of Theorems~\ref{thm3} and~\ref{thm-sheet}}

\subsection{Proof of Theorem~\ref{thm3}} \label{sec-proof-thm3}

We will use the results from~\cite{ChCh} and~\cite{FG1,FG2} (for $d=2$
and $d \geq 3$, respectively) cited in Section~\ref{sec-enh-dim}
and stated in the theorem below for convenience.
\begin{theorem}[\cite{ChCh},\cite{FG1,FG2}] \label{thm2}
For all $d \geq 2$, we have that $p_c(N,d) \geq p_c({\mathbb L}^d)$
for every $N \geq2$.
Moreover, $\lim_{N \rightarrow \infty}p_c(N,d)=p_c({\mathbb L}^d)$.
\end{theorem}

\noindent
With this, we are now ready to present the actual proof. \\

\noindent{\bf Proof of Theorem \ref{thm3}.} Consider Mandelbrot's
fractal percolation process with retention parameter $p$. Let
${\cal A}_k$ denote the event that there is complete retention
until level $k$ (as pointed out in~\cite{ChCh}, ``a minor miracle"),
i.e., $C_N^k=[0,1]^d.$ Trivially, for any $k\geq 1,$
\[
\theta_{N,d}(p)
\leq {\mathbb P}_p(CR(C_N) \, | \, {\cal A}_k).
\]

Assume that ${\cal A}_{k-1}$ has occurred. We will couple level $k$ of
the fractal percolation process, conditioned on ${\cal A}_{k-1}$, to
a diminishment percolation process on ${\mathbb L}^d$ with density of
diminishment $s=1-\theta_{N,d}(p)$. We now explain how the coupling works.

If ${\cal A}_{k-1}$ occurs, at the next level of the fractal percolation process,
the unit cube $[0,1]^d$ is divided in $N^{dk}$ subcubes of side length $1/N^k$.
We identify the subcubes with their centers
$\{ (\frac{x_1-1/2}{N^k},\ldots,\frac{x_d-1/2}{N^k}), 1 \leq x_1,\ldots,x_d \leq N^k \}$.
Each subcube is retained with probability $p$ and discarded with probability
$1-p$, independently of all other subcubes.

Recall that $\Lambda^d_N = \{1, \ldots, N \}^d \subset {\mathbb Z}^d$.
In our diminishment model, sites in ${\mathbb Z}^d \setminus \Lambda^d_{N^k}$
are declared closed, while a site $x=(x_1,\ldots,x_d) \in \Lambda^d_{N^k}$ is
declared open if the cube $(\frac{x_1-1/2}{N^k},\ldots,\frac{x_d-1/2}{N^k})$
is retained at level $k$ of the fractal percolation process, and closed
otherwise. This means that the states of sites in $\Lambda^d_{N^k}$
are determined by the fate of the level $k$ subcubes in the fractal
percolation process, while all other sites are closed.

On this percolation model we now perform a diminishment according
to the following rule. If a site
$x=(x_1,\ldots,x_d) \in \Lambda^d_{N^k}$ is open and all its
${\mathbb L}^d$-neighbors except possibly $(x_1 \pm 1,x_2,\ldots,x_d)$
are closed, we look at the effect of the fractal percolation process
within the subcube $(\frac{x_1-1/2}{N^k},\ldots,\frac{x_d-1/2}{N^k})$.
If the fractal percolation process inside that subcube does not leave a
connection between its two faces perpendicular to the first coordinate
direction, we declare site $x$ closed. The result is that, due to the
scale invariance of the fractal percolation process, each open site
$x=(x_1,\ldots,x_d) \in \Lambda^d_{N^k}$ such that all its
${\mathbb L}^d$-neighbors except possibly $(x_1 \pm 1,x_2,\ldots,x_d)$
are closed is independently declared closed with probability
$1-\theta_{N,d}(p)$. Sites that were already closed simply stay closed.
This is obviously a diminishment in the sense of Section~\ref{sec-enh-dim}.
We denote by ${\mathbb Q}^{dim}_{p,1-\theta_{N,d}(p)}$ the probability
measure corresponding to the diminishment percolation process just
defined (with initial configuration chosen as described above, i.e.,
with all sites outside $\Lambda^d_{N^k})$ closed).

We now make three simple but crucial observations. First of all, recalling
the definition of the diminishment percolation model of Section~\ref{sec-enh-dim}
and of the corresponding measure ${\mathbb P}^{dim}_{p,s}$, it is clear that,
due to the difference in the initial configurations before activating the
enhancement,
\begin{equation} \nonumber
{\mathbb Q}^{dim}_{p,s}(CR(\Lambda^d_{N^k},{\mathbb L}^d)) \leq
{\mathbb P}^{dim}_{p,s}(CR(\Lambda^d_{N^k},{\mathbb L}^d))
= \varphi^d_{N^k}(p,s)
\end{equation}
(where in our case $s = 1-\theta_{N,d}(p)$).

Secondly, note that any two retained level $k$ subcubes of the
fractal percolation process that touch at a single point (a corner)
will eventually not touch at all. The reason for this is that at
every subsequent step $j>k$ of the fractal construction, the corner
at which the two subcubes touch has probability uniformly bounded
away from zero to vanish because all the level $j$ subcubes that
share that corner are discarded. This means that, if we are
interested in the connectivity properties of $C_N$, any two level
$k$ subcubes that touch at a corner can effectively be considered
disconnected. (This is the reason why the lattice ${\mathbb L}^d$
is relevant in this context.)


Lastly, from topological considerations and the previous observation,
one can easily conclude that (conditioning on ${\cal A}_{k-1}$) the
absence of a crossing of $\Lambda({\mathbb L}^d,N^k)$ in the first
coordinate direction in the diminishment percolation process described
above implies the absence of a left to right crossing of $[0,1]^d$ by
the fractal set $C_N$. (Otherwise stated, if ${\cal A}_{k-1}$ happens,
then ``$CR(C_N) \subset CR(\Lambda^d_{N^k},{\mathbb L}^d)$," where
the event on the left concerns the fractal percolation process and
that on the right the diminishment percolation process, and the inclusion
should be interpreted in terms of the coupling between the two processes.)

From the previous discussion we can conclude that
\begin{eqnarray} \label{bound1-on-theta} \nonumber
\theta_{N,d}(p) & \leq & \P_p(CR(C_N) \, | \, {\cal A}_{k-1}) \\ \nonumber
& \leq & {\mathbb Q}^{dim}_{p,1-\theta_{N,d}(p)}(CR(\Lambda^d_{N^k},{\mathbb L}^d)) \\ \nonumber
& \leq & \P^{dim}_{p,1-\theta_{N,d}(p)}(CR(\Lambda^d_{N^k},{\mathbb L}^d)) \\
& = & \varphi^d_{N^k}(p,1-\theta_{N,d}(p)).
\end{eqnarray}

We will next use~(\ref{bound1-on-theta}) to show that, for any $\epsilon>0$,
there exists an $N_1=N_1(\epsilon)$ such that, for every $N \geq N_1$,
$\theta_N(p_c(N,d)+\tau) > 1-\epsilon$ for all $\tau>0$.
To do this, we set $s=\epsilon$ in Lemma~\ref{lemma-diminishment}, obtain
$\delta_2=\delta_2(\epsilon)>0$, and then use Theorem~\ref{thm2} to find
an $N_1=N_1(\delta_2)$ such that
\begin{equation} \label{eqn12} \nonumber
0 \leq p_c(N,d)-p_c({\mathbb L}^d)<\delta_2/2
\end{equation}
for all $N \geq N_1$.

Assume now, by contradiction, that
\begin{equation} \nonumber
\theta_{\hat N,d}(p_c(\hat N,d)+\hat\tau) \leq 1-\epsilon,
\end{equation}
for some $\hat N \geq N_1$ and some $\hat\tau$ which, because of
the monotonicity of $\theta_{N,d}(p)$ in $p$, we can assume without
loss of generality to be smaller than $\delta_2/2$.
Using the obvious monotonicity in $s$ of $\varphi^d_{N^k}(p,s)$,
from~(\ref{bound1-on-theta}) we obtain
\begin{equation} \label{bound2-on-theta}
\theta_{\hat N,d}(p_c(\hat N,d)+\hat\tau) \leq \varphi^d_{\hat{N}^k}(p_c(\hat N,d)+\hat\tau,\epsilon).
\end{equation}
Noting that, for any $\tau<\delta_2/2$,
\begin{equation} \nonumber
0 \leq p_c(\hat N,d)+\tau-p_c({\mathbb L}^d)<\delta_2,
\end{equation}
we can let $k \to \infty$ in~(\ref{bound2-on-theta}) and use
Lemma~\ref{lemma-diminishment} to obtain
\begin{equation} \label{theta=0}
\theta_{\hat N,d}(p_c(\hat N,d)+\hat\tau)=0.
\end{equation}

This however is an obvious contradiction, from which it follows
that for every $N \geq N_1$,
\begin{equation} \label{bound}
\theta_{N,d}(p_c(N,d)+\tau) > 1-\epsilon
\end{equation}
for all $\tau>0$, as claimed above.

We are now ready to conclude the proof. It is easy to see from the
construction of $C_N=\cap_{i=1}^{\infty}C_N^i$ that the function
$\theta_{N,d}(p)$, being the limit of decreasing continuous functions,
is right-continuous. Thus, letting $\tau \to 0$ in~(\ref{bound}),
we obtain
\begin{equation} \nonumber
\theta_{N,d}(p_c(N,d)) \geq 1-\epsilon
\end{equation}
for all $N \geq N_1$.
Since $\epsilon$ was arbitrary, this concludes the proof.
{QED}

\subsection{Proof of Theorem \ref{thm-sheet}} \label{sec-proof-thm-sheet}

We will use the following result from~\cite{Orzechowski} cited
in Section~\ref{sec-enh-dim} and stated in the theorem below
for convenience.
\begin{theorem}[\cite{Orzechowski}] \label{Orzechowski}
For all $d \geq 3$, we have that $\tilde{p}_c(N,d) \geq 1-p_c({\mathbb M}^d)$
for every $N \geq 2$.
Moreover, $\lim_{N \rightarrow \infty}\tilde{p}_c(N,d)=1-p_c({\mathbb M}^d)$.
\end{theorem}

\noindent
With this, we are now ready to present the actual proof. \\

\noindent{\bf Proof of Theorem \ref{thm-sheet}.}
Consider Mandelbrot's fractal percolation process with retention parameter $p$.
Let ${\cal A}_k$ denote the event that there is complete retention
until level $k$, i.e., $C_N^k=[0,1]^d.$ Trivially, for any $k\geq 1,$
\[
\theta_{N,d}(p)
\leq {\mathbb P}_p(CR(C_N) \, | \, {\cal A}_k).
\]

Assume that ${\cal A}_{k-1}$ has occurred. We will couple level $k$ of
the fractal percolation process, conditioned on ${\cal A}_{k-1}$,
to an enhancement percolation process on ${\mathbb M}^d$ with density
of enhancement $s=1-\tilde\theta_{N,d}(p)$ in a way similar to the
proof of Theorem~\ref{thm3}, as explained below.

If ${\cal A}_{k-1}$ occurs, at the next level of the fractal percolation process,
the unit cube $[0,1]^d$ is divided in $N^{dk}$ subcubes of side length $1/N^k$.
Like in the proof of Theorem~\ref{thm3}, we identify the subcubes with their centers
$\{ (\frac{x_1-1/2}{N^k},\ldots,\frac{x_d-1/2}{N^k}), 1 \leq x_1,\ldots,x_d \leq N^k \}$.
Each subcube is retained with probability $p$ and discarded with probability
$1-p$, independently of all other subcubes.

Contrary to the proof of Theorem~\ref{thm3}, in our enhancement model,
sites in ${\mathbb Z}^d \setminus \Lambda^d_{N^k}$ are declared open,
while a site $x=(x_1,\ldots,x_d) \in \Lambda^d_{N^k}$ is declared
closed if the cube $(\frac{x_1-1/2}{N^k},\ldots,\frac{x_d-1/2}{N^k})$ is
retained, and open otherwise. Once again, the states of sites in
$\Lambda^d_{N^k}$ are determined by the fate of the level $k$
subcubes of the fractal percolation process. This time however,
all sites outside $\Lambda^d_{N^k}$ are open.

On this percolation model we now perform an enhancement according to
the following rule. If a site $x=(x_1,\ldots,x_d) \in \Lambda^d_{N^k}$
is closed and $(x_1 + 1,x_2,\ldots, x_d)$ and $(x_1 - 1,x_2,\ldots, x_d)$
are both open, we look at the effect of the fractal percolation process
within the subcube $(\frac{x_1-1/2}{N^k},\ldots,\frac{x_d-1/2}{N^k})$.
If the fractal percolation process inside that subcube does not leave
a $(d-1)$-dimensional surface separating its two faces perpendicular to
the first coordinate direction, we declare site $x$ open. This time the
result is that each closed site $x=(x_1,\ldots,x_d) \in \Lambda^d_{N^k}$
such that its two neighbors $(x_1 + 1,x_2,\ldots, x_d)$ and
$(x_1 - 1,x_2,\ldots, x_d)$ are both open is independently declared open
with probability $1-\tilde\theta_{N,d}(p)$, while sites that were already
open simply stay open. This is obviously an enhancement in the sense of
Section~\ref{sec-enh-dim}. We denote by
${\mathbb Q}^{enh}_{1-p,1-\tilde\theta_{N,d}(p)}$ the probability measure
corresponding to the enhancement percolation process just defined (with
initial configuration chosen as described above, i.e., with all sites
outside $\Lambda^d_{N^k}$ open).

We now make two simple but crucial observations. First of all, recalling
the definition of the enhancement percolation model of Section~\ref{sec-enh-dim}
and of the corresponding measure ${\mathbb P}^{enh}_{p,s}$, it is clear that
\begin{equation} \nonumber
{\mathbb Q}^{enh}_{1-p,s}(CR(\Lambda^d_{N^k},{\mathbb M}^d)) \geq
{\mathbb P}^{enh}_{1-p,s}(CR(\Lambda^d_{N^k},{\mathbb M}^d))
= \psi^d_{N^k}(1-p,s)
\end{equation}
(where in this case $s = 1-\tilde\theta_{N,d}(p)$).

Secondly, from topological considerations one can easily conclude that
(conditioned on the event ${\cal A}_{k-1}$) the presence of an open
crossing of $\Lambda^d_{N^k}$ in the first coordinate in the
enhancement percolation process described above implies the absence of
a sheet crossing of $[0,1]^d$ by the fractal set $C_N$ separating the
two faces perpendicular to the first coordinate direction. (Otherwise
stated, if ${\cal A}_{k-1}$ happens, then
``$SH(C_N) \subset CR(\Lambda^d_{N^k},{\mathbb M}^d)^c$," where the event
on the left concerns the fractal percolation process and that on the right
the enhancement percolation process, the inclusion should be interpreted
in terms of the coupling between the two processes, and $CR(\cdot)^c$
denotes the complement of $CR(\cdot)$.)

From the previous discussion we can conclude that
\begin{eqnarray} \label{bound-on-theta-tilde} \nonumber
\tilde\theta_{N,d}(p) & \leq & \P_p(SH(C_N) \, | \, {\cal A}_{k-1}) \\ \nonumber
& \leq & 1 - {\mathbb Q}^{enh}_{1-p,1-\tilde\theta_{N,d}(p)}(CR(\Lambda^d_{N^k},{\mathbb M}^d)) \\ \nonumber
& \leq & 1 - \P^{enh}_{1-p,1-\tilde\theta_{N,d}(p)}(CR(\Lambda^d_{N^k},{\mathbb M}^d)) \\
& = & 1 - \psi^d_{N^k}(1-p,1-\tilde\theta_{N,d}(p)).
\end{eqnarray}

As in the proof of Theorem~\ref{thm3}, (\ref{bound-on-theta-tilde}) can be
used to show that, for any $\epsilon>0$, there exists an $N_2=N_2(\epsilon)$
such that, for every $N \geq N_2$,
$\tilde\theta_N(\tilde{p}_c(N,d)+\tau) > 1-\epsilon$ for all $\tau>0$.
To do this, we set $s=\epsilon$ in Lemma~\ref{lemma-enhancement}, obtain
$\delta_1=\delta_1(\epsilon)>0$, and then use Theorem~\ref{Orzechowski}
to find an $N_2=N_2(\delta_1(\epsilon))$ such that
\begin{equation} \nonumber
0 \leq \tilde{p}_c(N,d)-(1-p_c({\mathbb M}^d))<\delta_1/2,
\end{equation}
or equivalently,
\begin{equation} \nonumber
0 \leq p_c({\mathbb M}^d)-(1-\tilde{p}_c(N,d))<\delta_1/2,
\end{equation}
for all $N \geq N_2$.

From here, the proof can be concluded in essentially the same way as the
proof of Theorem~\ref{thm3}, with the difference that one needs to use
Lemma~\ref{lemma-enhancement} instead of Lemma~\ref{lemma-diminishment}.
We leave the details to the reader.
{QED}

\section{Bounds for $p_c(N,2).$} \label{sec-proof-2d}

In this section, we consider only two-dimensional fractal percolation,
and for simplicity we will not indicate the dimension in the notation.

Recall the definition of correlation length from Section~\ref{section2d}:
\[
L_{\delta}(p):=\min \{N: \varphi_N(p)\geq 1-\delta\},
\]
for any $p>p_c=p_c({\mathbb L}^2)=p_c({\mathbb Z}^2)$ and any $0<\delta<1/2$. \\

\noindent{\bf Proof of Theorem \ref{thm4}.} We start by proving (\ref{eqn15}).
One can easily show (see~\cite{ChCh}) that $\theta_N(p_c)=0$. Since $\theta_N(p_c(N))>0$
for every $N \geq 2$, we deduce that $p_c(N)>p_c$ for every $N \geq 2$.
Assume, by contradiction, that there exists a $C_0<\infty$ such that
\begin{equation} \label{eq-contra}
p_c(N)<p_c+C_0 \left( \frac{1}{N} \right)^{3/4}
\end{equation}
for all $N \geq 2$. Fix $\delta>0$; by assumption {\bf H1} and~(\ref{eq-contra}),
we have
$L_{\delta}(p_c(N)) \geq c_1(\delta) (p_c(N)-p_c)^{-4/3} > c_1(\delta) C_0^{-4/3} N$,
or equivalently, $N<c' L_{\delta}(p_c(N))$ for all $N \geq 2$ and some constant
$c'=c'(\delta)$. It then follows from Theorem~1 of~\cite{Kesten} or Theorem~26
of~\cite{Nolin} applied to the square lattice (see Section~8 of~\cite{Nolin}
concerning the applicability of such results to various regular lattices), together
with the Russo-Seymour-Welsh theorem, that $\varphi_N(p_c(N))$ is bounded away from
$1$ uniformly in $N$. However, that leads to a contradiction since
$\varphi_N(p_c(N)) \geq \theta_N(p_c(N))$ and $\lim_{N \to \infty}\theta_N(p_c(N))=1$
by Theorem~\ref{thm3}.

We now proceed to prove (\ref{eqn16}). First, let $R_m$ be a
$3m/N \times m/N$ rectangle and denote by ${\cal B}^1$ the event that
$R_m$ contains a left-to-right crossing by $C^1_N$, and both the left
and the right third of $R_m$ contain a top-to-bottom crossing by $C^1_N$
(see Appendix A for more motivation on why we consider this event, and
its relation with the functions $f_{N,m}$ and $g_{N,m}$ that will be
introduced below). For any $\delta_1>0$, one can choose $\delta_0>0$ so
small that the probability of the event ${\cal B}^1$ for any $p>p_c$ is
at least $1-\delta_1$ when $m$ is chosen to be $L_{\delta_0}(p)$.
To see this note that, letting $m=L_{\delta_0}(p)$, we have
$\varphi_m(p) \geq 1-\delta_0$. The claim then follows from standard
arguments using the Russo-Seymour-Welsh theorem and the Harris-FKG
inequality. Observe that the choice of $\delta_0$ is only dependent
on $\delta_1$.

For any $N\geq2$, let
\begin{equation} \label{eqn30}
p_N:=p_c+D\left(\frac{\log N}{N}\right)^{3/4},
\end{equation}
where $D$ is a constant yet to be determined.
Furthermore, we will let $m=L_{\delta_0}(p_N)$, for a $\delta_0>0$
to be chosen later. The definition of $p_N$ and hypothesis~{\bf H1}
imply that
\begin{equation} \label{bound-on-m}
m \leq \frac{c_2(\delta_0)}{D^{4/3}} \frac{N}{\log{N}}.
\end{equation}

Following~\cite{ChCh}, we now introduce the functions
\[
f_{N,m}(y) := \frac{c_3}{y (1 - g y^{1/4})} N^2 (g y^{1/4})^{N/2m}
\]
and
\[
g_{N,m}(y) := \frac{c_4}{y (1 - g y^{1/4})} \frac{N}{m} (g y^{1/4})^{N/2m}
= \frac{c_4}{c_3} \frac{1}{Nm} f_{N,m}(y)
\]
where $c_3,c_4>0$ are finite constants and $g \leq 3$.
The role of these two functions will be clear soon. For completeness and
the reader's convenience, their derivation is presented in Appendix A.

Fix a $y_0>0$ such that $g y_0^{1/4} < 1$, let $\delta_1=y_0/2$, and choose
$\delta_0>0$ so that $\P_p({\cal B}^1) \geq 1 - \delta_1$ for any $p>p_c$.
Choose $D$ in the definition of $p_N$ large enough so that
\[
h_N(y):=\frac{c_5}{y}N^2(g y^{1/4})^{\frac{D^{4/3}}{2 c_2(\delta_0)} \log N}
\rightarrow 0, \ \ \forall y \leq y_0 \ \ \textrm{ as } N\rightarrow\infty,
\]
where $c_5 = \frac{c_3}{1 - g y_0^{1/4}}$. Recalling that $m=L_{\delta_0}(p_N)$,
we can then use~(\ref{bound-on-m}) to conclude that, for all $y \leq y_0$,
\begin{equation} \label{eqn-f-h}
f_{N,m}(y) \leq h_N(y).
\end{equation}

With the choices made above, it is easy to see that, for $N$ sufficiently large,
the equation
\begin{equation} \label{h-fixed-point}
y=\delta_1+h_N(y)
\end{equation}
has a solution $y_2(N)$ which converges to $\delta_1$ as $N \to \infty$.
Moreover, (\ref{eqn-f-h}) implies that, whenever~(\ref{h-fixed-point}) has
a solution $y_2(N)$,
\begin{equation} \label{f-fixed-point}
y=\delta_1+f_{N,m}(y)
\end{equation}
has a solution $y_1(N)$ such that $y_1(N) \leq y_2(N)$.
We pick $N_1$ such that, for all $N \geq N_1$, (\ref{h-fixed-point}) has a
solution and $h_N(y)$ is increasing in $y$.

We can now use a bound from~\cite{ChCh}, rederived in Appendix A
(see~(\ref{lower-bound}) there), to show that
\begin{eqnarray} \nonumber
\P_{p_N}(CR(C_N^k)) & \geq & 1 - g_{N,m}(\delta_k) \\
& = & 1 - \frac{c_4}{c_3} \frac{1}{Nm} f_{N,m}(\delta_k) \nonumber \\
& \geq & 1 - \frac{c_4}{c_3}\frac{h_{N}(\delta_k)}{N}, \label{array-lower-bound}
\end{eqnarray}
where $\{ \delta_k \}_{k \geq 2}$ is some sequence that satisfies
the iterative inequality
\begin{equation} \label{iter-eq}
\delta_{k} \leq \delta_1 + f_{N,m}(\delta_{k-1}).
\end{equation}
In order for~(\ref{array-lower-bound}) to hold, it suffices that
$\P_{p_N}({\cal B}^1) \geq 1 - \delta_1$, which is satisfied for
our choice of $p_N$ and $m(=L_{\delta_0}(p_N))$.
When~(\ref{f-fixed-point}) has a solution $y_1$, the sequence
$\{ \delta_k \}_{k \in {\mathbb N}}$ must have a limit $\delta^* \leq y_1$.

For $N \geq N_1$, letting $k \to \infty$ on both sides of~(\ref{array-lower-bound}),
and using~(\ref{h-fixed-point}) and the fact that $\delta^* \leq y_1 \leq y_2$
and $h_N(y)$ is increasing in $y$, we obtain
\begin{eqnarray}
\theta_N(p_N) & \geq & 1-\frac{c_4}{c_3}\frac{h_{N}(\delta^*)}{N} \nonumber \\
& \geq & 1-\frac{c_4}{c_3}\frac{h_{N}(y_2)}{N} \nonumber \\
& = & 1-\frac{c_4}{c_3}\frac{(y_2-\delta_1)}{N}. \nonumber
\end{eqnarray}

Since $\lim_{N \to \infty}y_2(N)=\delta_1$, as remarked above, we can conclude
that there exists an $N_1$ (possibly larger than required before) such that
\[
\theta_N(p_N) \geq 1-\frac{c_4}{c_3}\frac{(y_2-\delta_1)}{N} > 0
\]
for all $N \geq N_1$, which concludes the proof.
{QED} \\


\noindent {\bf Acknowledgements.} It is a pleasure to thank Ronald Meester for many
interesting discussions on subjects closely related to the content of this paper.
The interest of one of the authors (F.C.) in Mandelbrot's fractal percolation process
was started by a question of Alberto Gandolfi during the workshop ``Stochastic Processes
in Mathematical Physics" at Villa La Pietra, Florence, June 19-23, 2006.
The other author (E.B.) thanks everyone at the Vrije Universiteit Amsterdam and Ronald Meester
in particular for making his one year stay possible as well as enjoyable and productive.

\appendix
\refstepcounter{section}
\section*{Appendix A} \label{ChCh}

For completeness and for the reader's convenience, following~\cite{ChCh},
in this appendix we present the derivation of the functions $f_{N,m}$ and
$g_{N,m}$ used in the proof of Theorem~\ref{thm4} (note that $g_{N,m}$ is
not given a name in~\cite{ChCh}).

Let $m$ be a fixed number much smaller than $N$, and denote by $H_m(x_1,x_2)$
the $3m/N \times m/N$ rectangle
$\{ y=(y_1,y_2) \in {\mathbb R}^2 : x_1 \leq y_1 \leq x_1+3m/N, x_2 \leq y_2 \leq x_2+m/N \}$,
and by $V_m(x_1,x_2)$ the $m/N \times 3m/N$ rectangle
$\{ y=(y_1,y_2) \in {\mathbb R}^2 : x_1 \leq y_1 \leq x_1+m/N, x_2 \leq y_2 \leq x_2+3m/N \}$.
If at level $k$ of the fractal percolation process $H_m(x_1,x_2)$ contains
a left-to-right crossing by $C^k_N$, and both the left and the right third
of $H_m(x_1,x_2)$ contain a top-to-bottom crossing by $C^k_N$, we say that
the event ${\cal H}^k(x_1,x_2)$ happens. Analogously, if at level $k$
of the fractal percolation process $V_m(x_1,x_2)$ contains a top-to-bottom
crossing by $C^k_N$, and both the top and the bottom third of $V_m(x_1,x_2)$
contain a left-to-right crossing by $C^k_N$, we say that ${\cal V}^k(x_1,x_2)$
happens.

\begin{figure}[!ht]
\begin{center}
\includegraphics[width=5cm]{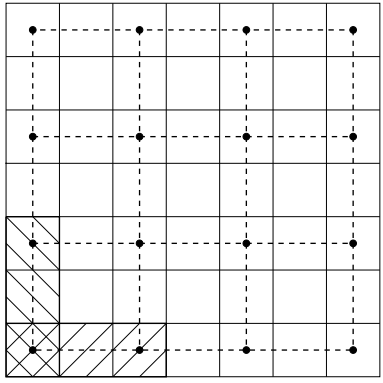}
\caption{The unit square is partially covered by overlapping $3m/N \times m/N$ and $m/N \times 3m/N$
rectangles like the shaded ones shown in the bottom-left corner. The centers of the squares
where the rectangles overlap, indicated by dots, are the sites of an associated portion of
the square lattice whose bonds are indicated by broken segments.}
\label{fig1}
\end{center}
\end{figure}

We can partially cover the unit square $[0,1]^2$ using rectangles $H_m(x_1,x_2)$
and $V_m(x_1,x_2)$ with $x_1$ taking values $0, 2m/N, 4m/N, \ldots$ and $x_2$ taking
values $0, 2m/N, 4m/N, \ldots$ . We call $\cal R$ the set of overlapping rectangles
used for this partial covering of the unit square, and assume for simplicity that
$N/m$ is an odd integer (see Figure~\ref{fig1}). The readers can easily convince
themselves that this assumption can be safely removed in the arguments below.
The probabilities of the crossing events ${\cal H}^k(\cdot,\cdot)$ and
${\cal V}^k(\cdot,\cdot)$ are equal and are the same for all rectangles in
$\cal R$. We therefore denote them $\pi_k$.

We identify the $m/N \times m/N$ squares where the rectangles from $\cal R$
overlap with the sites of a portion of the square lattice (see Fig.~\ref{fig1}).
On this portion of the square lattice, we consider a family of bond percolation
models indexed by $k=1,2,\ldots$ . In the $k$-th member of the family,
we declare a horizontal bond between two adjacent sites open if the event
${\cal H}^k(x_1,x_2)$ happens, where $(x_1,x_2)$ is the bottom-left corner of the
$m/N \times m/N$ square corresponding to the leftmost of the two adjacent sites.
Analogously, we declare a vertical bond between two adjacent sites open if the event
${\cal V}^k(x_1,x_2)$ happens, where $(x_1,x_2)$ is the bottom-left corner of the
$m/N \times m/N$ square corresponding to the lower site. Note that we obtain a family
of dependent bond percolation models with densities $\pi_k, k=1,2,\ldots$ , of open
bonds, and that a left-to-right crossing in the $k$-th bond percolation process implies
a left-to-right crossing of the unit square by $C^k_N$. Note also that for $k=1$, the
event ${\cal H}^1(\cdot,\cdot)$ (resp., ${\cal V}^1(\cdot,\cdot)$) is simply a crossing
event for Bernoulli site percolation on a $3m \times m$ (resp., a $m \times 3m$) rectangular
portion of the square lattice.

\begin{figure}[!ht]
\begin{center}
\includegraphics[width=9cm]{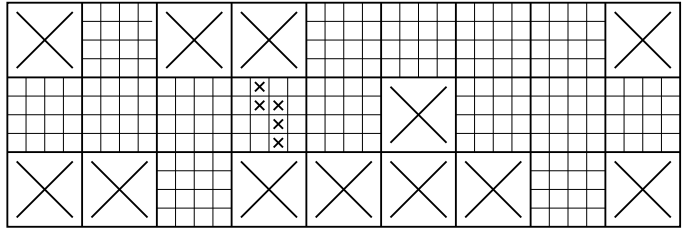}
\caption{A $3m/N \times m/N$ rectangle where the event ${\cal B}^1$ happens but ${\cal B}^2$
does not. The larger squares are of size $1/N \times 1/N$. Squares with a cross are not
retained. Retained squares are divided in smaller squares of size $1/N^2 \times 1/N^2$.
Here, $m=3$ and $N=4$.}
\label{fig2}
\end{center}
\end{figure}

Let us now consider for definiteness a given rectangle $R_m$ from $\cal R$. We denote
by ${\cal B}^k$ the relevant crossing event (either ${\cal H}^k(\cdot,\cdot)$ or
${\cal V}^k(\cdot,\cdot)$) inside $R_m$ at the $k$:th step of the fractal percolation
process. We have $\P_p({\cal B}^k)=\pi_k$. Let us assume that
$\pi_1 \geq 1 - \delta_1$ for some fixed $\delta_1$. Assuming that ${\cal B}^1$ happens,
we are going to obtain a lower bound for $\P_p({\cal B}^2)$ in terms of $\delta_1$
(see Fig.~\ref{fig2}).

After the second step of the fractal percolation process, $R_m$ is split into $3m^2 \times N^2$
squares of side length $1/N^2$. We can partially cover $R_m$ using rectangles $H_{m/N}(\cdot,\cdot)$ and
$V_{m/N}(\cdot,\cdot)$ in the same fashion as we did for the unit square. Each such rectangle
contains $3m^2$ squares of side length $1/N^2$. If $R_{m/N}$ is a rectangle contained inside a
$1/N \times 1/N$ square that was retained in the first step of the fractal percolation process,
the probability to see inside $R_{m/N}$ a scaled copy of the event ${\cal B}^1$ is exactly
$\pi_1$ because of the scale invariance of the fractal percolation process. More generally,
as we did above inside the unit square, we can couple step $k+1$ of the fractal percolation
process inside the retained portion of $R_m$ to a dependent bond percolation process with
density of open bonds $\pi_k$.

We will now use the bond percolation process coupled to the second level of the fractal
percolation process to obtain a lower bound for $\P_p({\cal B}^2 \, | \, {\cal B}^1)$ in
terms of $\delta_1$. One way to ensure the occurrence of ${\cal B}^2$ is by a string of
open bonds in the coupled bond percolation process. The event that such a string is not
present is by duality equivalent to the presence of a string of closed dual bonds, where
a dual bond is perpendicular to an ordinary bond and is declared closed (resp., open) if
the (unique) ordinary bond it intersects is closed (resp., open). By simple
worst-case-scenario counting arguments, one can give an upper bound for the probability
of such a string of closed dual bonds, as explained below. Note that we are considering
a dependent bond percolation model, but that the dependence range is finite (and in fact
very short). We can certainly assume, as in~\cite{ChCh}, that finding a string consisting
of $l$ closed dual bonds implies at least $l/4$ independent trials.

By using trivial upper bounds for the number of strings of length $l$, we obtain
\begin{eqnarray} \nonumber
\P_p({\cal B}^2 \, | \, {\cal B}^1) & \geq & 1 - 3m^2 \left(\frac{4N}{m}\right)
\sum_{l \geq \frac{N}{2m}-1} g^l \delta_1^{l/4} \\
& \geq & 1 - c'_3 N m \frac{(g \delta_1^{1/4})^{N/2m}}{g \delta_1^{1/4} (1 - g \delta_1^{1/4})} \nonumber \\
& \geq & 1 - c_3 N^2 \frac{(g \delta_1^{1/4})^{N/2m}}{\delta_1 (1 - g \delta_1^{1/4})} \nonumber
\end{eqnarray}
where $g \leq 3$ is an upper bound on the number of possible choices for attaching the
next bond to a growing string, and $c'_3,c_3>0$.

Note that in the first line of the above equation, the factor $\frac{4N}{m}$ comes from
the fact that the string of closed dual bonds must start at one of the four sides of a
retained $1/N \times 1/N$ square, while $3m^2$ is a trivial upper bound on the number of
retained $1/N \times 1/N$ squares, which can certainly be no bigger than the total number
of $1/N \times 1/N$ squares in $R_m$ (see Fig.~\ref{fig2}). The summation is over values
of $l$ large enough so that a string of length $l$ can prevent the event ${\cal B}^2$.
Introducing the function $f_{N,m}(y) := \frac{c_3}{y (1 - g y^{1/4})} N^2 (g y^{1/4})^{N/2m}$,
we can write $\P_p({\cal B}^2 \, | \, {\cal B}^1) \geq 1 - f_{N,m}(\delta_1)$.

Defining $\delta_2$ via the equation $\pi_2=\P_p({\cal B}^2)=1-\delta_2$, from the
discussion above we have
\begin{equation} \nonumber
1-\delta_2 = \P_p({\cal B}^2 \, | \, {\cal B}^1) \P_p({\cal B}^1)
\geq (1 - f_{N,m}(\delta_1)) (1 - \delta_1) \geq 1- \delta_1 - f_{N,m}(\delta_1).
\end{equation}
Next, we are going to obtain a lower bound for $\P_p({\cal B}^3 \, | \, {\cal B}^1)$
in terms of $\delta_2$.

One way to ensure the occurrence of ${\cal B}^3$ is by a string of open bonds in the
dependent bond percolation process coupled to the third step of the fractal percolation
process inside the portion of $R_m$ retained after the first step of the fractal percolation
process. The density of open bonds is now $\pi_2$. We can then repeat exactly the same arguments
as above, but with $\pi_2$ instead of $\pi_1$ as density of open bonds. This obviously gives
us $\P_p({\cal B}^3 \, | \, {\cal B}^1) \geq 1 - f_{N,m}(\delta_2)$ and
$1-\delta_{3} \geq 1- \delta_1 - f_{N,m}(\delta_2)$, where $\delta_3$ is defined via the
equation $\pi_3=\P_p({\cal B}^3)=1-\delta_3$.
Proceeding by induction, we obtain the iterative inequality
$1-\delta_{k+1} \geq 1- \delta_1 - f_{N,m}(\delta_k)$, i.e.,
\begin{equation} \nonumber
\delta_{k+1} \leq \delta_1 + f_{N,m}(\delta_k).
\end{equation}

A way to ensure a left-to-right crossing of the unit square by $C^k_N$ is by a left-to-right
crossing of open bonds in the dependent bond percolation process coupled to step $k$ of the
fractal percolation process inside the unit square. Therefore, using arguments analogous to
those used above, one can obtain a lower bound for the probability of a left-to-right crossing
of the unit square by $C^k_N$ in terms of $\delta_k$ (or $\pi_k=1-\delta_k$). Counting arguments
and approximations analogous to those used before give
\begin{equation} \label{lower-bound}
\P_p(CR(C_N^k)) \geq 1 - g_{N,m}(\delta_k),
\end{equation}
where
$g_{N,m}(y) := \frac{c_4}{y (1 - g y^{1/4})} \frac{N}{m} (g y^{1/4})^{N/2m}$, with $c_4$
a positive constant.


\bigskip


\begin{thebibliography}{99}

\bibitem{AG} M.~Aizenmann and G.~Grimmett, Strict Monotonicity for
Critical Points in Percolation and Ferromagnetic Models,
\emph{J.~Stat.~Phys.} {\bf 63}, 817--835 (1991).

\bibitem{bcks} C.~Borgs, J.~Chayes, H.~Kesten, J.~Spencer,
The Birth of the Infinite Cluster: Finite-Size Scaling in Percolation,
\emph{Comm.~Math.~Phys.} {\bf 224}, 153--204 (2001).

\bibitem{cfn1} F.~Camia, L.~R.~Fontes, C.~M.~Newman,
The Scaling Limit Geometry of Near-Critical 2D Percolation,
\emph{J.~Stat.~Phys.} {\bf 125}, 1155-1171 (2006).

\bibitem{cfn2} F.~Camia, L.~R.~Fontes, C.~M.~Newman,
Two-Dimensional Scaling Limits via Marked Nonsimple Loops,
\emph{Bull.~Braz.~Math.~Soc.} {\bf 37}, 537-559 (2006).

\bibitem{lchayes} L.~Chayes, Aspects of the fractal percolation process,
\emph{Progress in Probability} {\bf 37}, 113--143 (1995).

\bibitem{ChCh} J.T.~Chayes and L.~Chayes, The large-N limit of the
threshold values in Mandelbrot's fractal percolation process,
\emph{J.~Phys.~A: Math.~Gen.} {\bf 22}, L501--L506 (1989).

\bibitem{CCD} J.T.~Chayes, L.~Chayes and R.~Durrett, Connectivity
Properties of Mandelbrot's Percolation Process,
\emph{Probab.~Theory~Relat.~Fields} {\bf 77}, 307--324 (1988).

\bibitem{CCGS} J.T.~Chayes, L.~Chayes, E.~Grannan and G.~Swindle,
Phase transitions in Mandelbrot's percolation process in three
dimensions, \emph{Probab.~Theory~Relat.~Fields} {\bf 90}, 291--300
(1991).

\bibitem{ChN} L.~Chayes and P.~Nolin, Large Scale Properties of
the IIIC for 2D Percolation, preprint arXiv:0705.357 (2007).

\bibitem{DG} F.M.~Dekking and G.R.~Grimmett, Superbranching processes
and projections of random Cantor sets, \emph{Probab.~Theory~Relat.~Fields}
{\bf 78}, 335--355 (1988).

\bibitem{DM} F.M.~Dekking and R.W.J.~Meester, On the structure of
Mandelbrot's percolation process and other Random Cantor sets,
\emph{J.~Stat.~Phys.} {\bf58}, 1109--1126 (1990).

\bibitem{FG1} K.J.~Falconer and G.R.~Grimmett, The critical point of
fractal percolation in three and more dimensions,
\emph{J.~Phys.~A:~Math.~Gen.} {\bf 24}, L491--L494 (1991).

\bibitem{FG2} K.J.~Falconer and G.R.~Grimmett, On the geometry of Random Cantor
Sets and Fractal Percolation, \emph{J.~Theor.~Probab.} {\bf 5}, 465--485 (1992).

\bibitem{G} G.~Grimmett, {\em Percolation},
Second edition, Springer-Verlag, Berlin (1999).

\bibitem{Kesten} H.~Kesten, Scaling Relations for 2D-Percolation,
\emph{Commun.~Math.~Phys.} {\bf 109}, 109--156 (1987).

\bibitem{Mandelbrot} B.B.~Mandelbrot, \emph{The Fractal Geometry of
Nature}, W.H.~Freeman, San Francisco (1983).

\bibitem{Meester} R.W.J.~Meester, Connectivity in fractal percolation,
\emph{J.~Theor.~Prob.} {\bf 5}, 775--789 (1992).


\bibitem{MPV} M.V.~Menshikov, S.Yu.~Popov and M.~Vachkovskaia,
On the connectivity properties of the complementary set in fractal percolation models,
\emph{Probab.~Theory~Relat.~Fields} {\bf 119}, 176--186 (2001).

\bibitem{Nolin} P.~Nolin, Near-critical percolation in two dimensions, preprint
arXiv:0711.4948.

\bibitem{nw} P.~Nolin, W.~Werner, Asymmetry of near-critical percolation interfaces,
preprint arXiv:0710.1470.

\bibitem{Orzechowski} M.E.~Orzechowski, On the Phase Transition to
Sheet Percolation in Random Cantor Sets, \emph{J.~Stat.~Phys.} {\bf
82}, 1081--1098 (1996).

\bibitem{SW} S.~Smirnov and W.~Werner, Critical exponents for two-dimensional percolation,
\emph{Math.~Res.~Lett.} {\bf 8}, 729--744 (2001).

\bibitem{White} D.G.~White, On fractal percolation in ${\mathbb R}^2$,
 \emph{Statist. Probab. Lett.} {\bf 45}, 187--190 (1999).

\end{thebibliography}
\end{document}